%% file: draft1.tex
\documentclass[11pt]{article}
\usepackage{amssymb}
\usepackage{float}
\usepackage{amsmath}
\usepackage{graphicx}
\usepackage{amsthm}
\usepackage{mathtools}
\newtheorem{proposition}{Proposition}
\newtheorem{lemma}{Lemma}
\newtheorem{remark}{Remark}

\def\R{{I \kern -.30em R}}
\def\C{{I \kern -.57em C}}

\title{Perfectly matched layers in time domain. A simple two-dimensional error analysis.}

\author{
K.Bryan
\thanks{Terre Haute, Indiana ({\tt bryan@rose-hulman.edu}).}
 \and
M.S. Vogelius
\thanks{Department of Mathematics, Rutgers University,
New Brunswick, NJ 08903, USA,
({\tt vogelius@math.rutgers.edu}).} }

\begin{document}
\bibliographystyle{plain}
\baselineskip 15pt
\vskip -150pt

\maketitle

%%%%%%%%%%%%%%%%%%%%%%%%%%%%%%%%%%%%%%%%%%%%%%%%%%%%%%%%%%%%%%%%%%%%

\centerline{Dedicated to the memory of Ivo Babuska and his seminal}
\centerline{contributions to Applied and Numerical Analysis.}

\begin{abstract}
Perfectly Matched Layers (PML) has become a very common method for the numerical approximation
of wave and wave-like equations on unbounded domains. This technique allows one to obtain accurate
solutions while working on a finite computational domain, and the technique is relatively simple to implement.
Results concerning the accuracy of the PML method have been obtained, but mostly with regard problems
at a fixed frequency. In this paper we provide very explicit time-domain bounds on the accuracy of PML for
the two-dimensional wave equation, and illustrate our conclusions with some numerical examples.
\end{abstract}

%%%%%%%%%%%%%%%%%%%%%%%%%%%%%%%%%%%%%%%%%%%%%%%%%%%%%%%%%%%%%%%%%%%%

%%%%%%%%%%%%%%%%%%%%%%%%%%%%%%%%%%%%%%%%%%%%%%%%%%%%%%%%%%%%

\pagestyle{myheadings}
\mark{{\mbox{ }}Bryan and Vogelius \quad PML}

%%%%%%%%%%%%%%%%%%%%%%%%%%%%%%%%%%%%%%%%%%%%%%%%%%%%%%%%%%%%

\input{Intro1.tex}

\input{Inhom-wave-eq1.tex}
\input{numerics.tex}
\input{Appendix1}
%\input{Comp.tex}
%\input{ny3.3.tex}
%\input{ny4.8.tex}
%\input{ny5.4.tex}
%\input{ny7.3.tex}

\end{document}

%% file: Intro1.tex
%%%%%%%%%%%%%%%%%%%%%%%%%%% main body %%%%%%%%%%%%%%%%%%%%%%%%%%%%%%%%%
\section{Introduction}

Perfectly Matched Layers (PML) is a technique first proposed in \cite{Berenger1, Berenger2} to address the challenge of numerically
approximating solutions to Maxwell's equations and related PDEs posed on unbounded domains, but in which the solution
is of interest only on some bounded subdomain. The challenge is to limit computation to this smaller subdomain or
some modestly larger region. For example, we may seek a solution to the wave equation with compactly supported initial data on $\mathbb{R}^n$, but the only region of interest for the solution is a ball $B$ centered at the origin. This might
be approached by enclosing $B$ in a larger ball $B_R$ of radius $R$ and then solving the wave equation on $B_R$ with some
boundary conditions, for example, zero Dirichlet conditions. However for sufficiently large times this results in reflections at the computational domain boundary $\partial B_R$ that eventually corrupt the full-space solution in $B$. Attempting to circumvent this by making $R$ larger may result in a great deal of unnecessary work. PML provides a method for altering and extending the PDE beyond the region of interest in such a way that solutions in $B$ are unchanged but attenuate quickly
in the ``PML region'' outside $B$; moreover, no reflections are
generated at the interface between $B$ and the PML region. Since solutions are quickly attenuated in the PML region, the process
of solving the modified PDEs on a truncated domain $B_R$ results in minimal alteration of the solution in $B$. One can thus
produce accurate solutions to the full-space PDE while working on a finite computational domain.

PML has gained popularity because it is relatively straightforward to implement, efficient, and in practice performs well.
It has been applied in numerous contexts, e.g., the acoustic wave equation \cite{grotesim,kaltenbacher}, Maxwell's equations \cite{becache,Berenger1,Berenger2}, seismic wave analysis
\cite{komatrump}, the linearized Euler's equation in fluids \cite{hu}, and other areas. The implementation of PML for the monochromatic reduced wave equation, the Helmholtz equation, has certain similarities to the passive ``cloaking by mapping technique'' employed to achieve approximate invisibility.

The application of the PML technique to a scalar PDE or a system of PDEs typically results in the addition of
some number of auxiliary functions, as well as additional equations. A variety of interesting questions arise, for example, the issue of whether the resulting system well-posed and stable; see for example \cite{becache} or \cite{grotesim}.
The numerical methods used to solve the PML PDEs can also have an effect on the performance of this method, since the
reflectionless property enjoyed at the PML interface in the continuous case may not hold perfectly when discretized.
See \cite{durukreiss} for a survey of progress on this question.

Of paramount importance for computation is the issue of how well the solutions to the PML-modified system (or a version posed on a finite computational domain) approximate the solution to the original PDE on the unbounded domain.
Analysis of the accuracy of PML for approximation of the full-space solution to specific PDEs has been addressed and some rigorous error estimates have been derived, but primarily for the monochromatic (single frequency) reduced wave equation (the Helmholtz equation) and other fixed frequency cases; see for example \cite{lassas,Monk}.
To the best of our knowledge the situation is somewhat different in the time-domain. There have been several numerical studies for the wave equation \cite{kaltenbacher, Monk}, but few explicit analytical estimates have appeared. One notable result
is \cite{diazjoly}, in which the authors provide a time-domain estimate for the accuracy of the PML technique applied
to the wave equation when the region of interest is the left half-plane and a PML region is introduced in the right half-plane.
Their Theorem 4 provides a time-domain estimate with regard to the supremum norm for agreement between the full-space wave equation and the PML modified system.

Our contribution in this paper is to provide a very simple derivation of error estimates in the time domain when the region
of interest is the unit disk in the plane, based on Fourier analysis, the construction of explicit solutions, and straightforward estimates for the monochromatic case.

Specifically, as our time-dependent problem we take a $2+1$ dimensional problem of the form
$u_{tt}-\Delta u = f(t,x)$ in $(0,\infty)\times \R^2$, the inhomogeneous wave equation, with initial data
$u(0,x)=u_t(0,x)=0$ for $x\in\R^2$. We assume that for all $t$, the function $f(t,\cdot)$ is compactly supported in the unit ball $B_1$ centered
at the origin.
The goal is to approximate the solution $u(t,x)$ for $x\in B_1$ while performing computations only on a ball
of radius $R>1$.

In Section \ref{sec:waveplane} we work in $(r,\theta)$ polar coordinates to construct explicit solutions to this
inhomogeneous wave equation, and then in
Section \ref{sec:pmlsolution} we consider a PML-modified version of these solutions. These modified solutions
are defined on a ball of radius $R>1$ and have zero Dirichlet data
on $r=R$. Moreover, these solutions agree closely and quantifiably in the supremum norm
with the full-space solution $u(t,x)$ for $x$ in the unit ball for $t>0$.
These solutions are constructed using the ``complex change-of-coordinates'' approach (see \cite{johnson} for a simple
treatment of this PML technique). In Section \ref{sec:example} we provide an example and general conclusion concerning
the efficacy of PML in the time domain. In Sections \ref{sec:pmlpdes} and \ref{sec:initbound} we show that these PML-modified solutions satisfy a certain system of PDEs in the plane with specific initial data. Finally in Section \ref{sec:numerics}
we illustrate our conclusions with a few simple numerical experiments.

%% file: Inhom-wave-eq1.tex
%%%%%%%%%%%%%%%%%%%%%%%%%%% main body %%%%%%%%%%%%%%%%%%%%%%%%%%%%%%%%%
\section{The Inhomogeneous wave equation on the plane}
\label{sec:waveplane}

In this section we shall first derive an explicit expression for the solution to the following initial value problem for the linear wave equation with wave speed $c=1$,
\begin{eqnarray}
\label{wave-eq}
u_{tt}-\Delta u = f(t,x)~ \hbox{ in }~ (0,\infty)\times \R^2~, \\
u(0,x)=u_t(0,x)=0 ~\hbox{ for } x\in \R^2~. \nonumber
\end{eqnarray}
We shall take $f(t,x)$ to be supported in $(0,\infty)\times \{\, |x|\le \frac12 \,\}$ and of the special form
\begin{equation}
\label{eqn:fdef}
f(t,r,\theta)= \sum_{n=0}^N \cos n\theta \, f_n(t,r)
\end{equation}
where $(r,\theta)$ are standard polar coordinates.
If we extend both $u$ and $f$ to be zero for $t\le 0$  ({\it i.e.}, all the $f_n$ vanish for $t\le 0$) then $u$ solves (\ref{wave-eq}) in $(-\infty,\infty)\times \R^2$.

We define functions $\hat f(\omega,x)=\frac{1}{\sqrt{2\pi}}\int_0^\infty e^{i\omega t}f(t,x)~dt$ and $\hat u(\omega,x)=\frac{1}{\sqrt{2\pi}}\int_0^\infty e^{i\omega t}u(t,x)~dt$, $\omega \in \R$. We assume that $f$ is real-valued, and therefore $u$ as well, and so we have
$$
\hat u(-\omega,x)= \overline{\hat u(\omega,x)}~~\hbox{ and } ~~ \hat f(-\omega,x)= \overline{\hat f(\omega,x)}~.
$$
The standard inversion formula for the Fourier Transform (in the case of $u$) may now be rewritten as
$$
u(t,x) = 2 \Re \left(\frac{1}{\sqrt{2\pi}} \int_0^\infty e^{-i\omega t} \hat u(\omega,x)\, d\omega   \right)~.
$$
For these reasons it suffices to only consider $\hat u(\omega,x)$ for $\omega>0$.
From (\ref{wave-eq}) it follows that the functions $\hat u(\omega,x)$ and $\hat f(\omega,x)$, $\omega>0$, are related by
$$
-\omega^2 \hat u(\omega,x)-\Delta_x \hat u(\omega,x) = \hat f(\omega,x) ~~ x\in \R^2~,
$$
with $\hat u$ satisfying the outgoing radiation condition (see \cite{Nguyen-Vogelius} for a derivation of this).
From $\Delta_x=\frac{1}{r}\left (\frac{\partial}{\partial r} \left (r\frac{\partial}{\partial r}\right )\right)
-\frac{1}{r^2}\frac{\partial^2}{\partial\theta^2}$ and the form of $f$ in (\ref{eqn:fdef}) we conclude that
\begin{equation}
\label{uhatsum}
\hat u(\omega,x)=\hat u(\omega,r,\theta)= - \sum_{n=0}^N \cos n \theta \, \hat \phi_n(\omega,r)~,
\end{equation}
where each function $\hat \phi_n(\omega,r)$ satisfies
$$
\frac{1}{r}\frac{\partial}{\partial r}\left( r\frac{\partial}{\partial r}\hat \phi_n\right)+\left(\omega^2-\frac{n^2}{r^2}\right)\hat \phi_n= \hat f_n(\omega,r)~,
$$
in addition to the outgoing radiation condition.

We define function $\tilde \psi_n$ by the relationship $\hat \phi_n(\omega,r)=\tilde \psi_n(\omega,\omega r)$.
The function $\tilde \psi_n$ satisfies the inhomogeneous Bessel equation
$$
\frac{1}{r}\frac{\partial}{\partial r}\left( r\frac{\partial}{\partial r}\tilde \psi_n\right)+\left(1-\frac{n^2}{r^2}\right)\hat \psi_n= \frac{1}{\omega^2}\hat f_n(\omega,\frac{r}{\omega})~,
$$
or
\begin{equation}
\label{eqn:besselinhom}
\frac{\partial^2}{\partial r^2} \tilde \psi_n +\frac{1}{r}\frac{\partial}{\partial r} \tilde \psi_n +\left(1-\frac{n^2}{r^2}\right)\tilde \psi_n = \frac{1}{\omega^2} \hat f_n(\omega,\frac{r}{\omega})~,
\end{equation}
as well as the outgoing radiation condition (with $\omega=1$). Note that the Bessel functions $J_n(r)$ and $Y_n(r)$
provide independent solutions to the homogeneous version of (\ref{eqn:besselinhom}), with Wronskian $J_n(r)Y'_n(r)-J'_n(r)Y_n(r)=2/(\pi r)$.
A straightforward variation of parameters computation then shows that solutions to (\ref{eqn:besselinhom}) take the form
\begin{eqnarray*}
\tilde \psi_n(\omega,r) &=&\frac{\pi}{2}Y_n(r)\int_0^r \frac{s J_n(s)}{\omega^2}\hat f_n(\omega,\frac{s}{\omega})\,ds -\frac{\pi}{2} J_n(r)\int_0^{r}\frac{s Y_n(s)}{\omega^2}\hat f_n(\omega,\frac{s}{\omega})\,ds \\
&&\hskip 75pt +c_{1,n,\omega} J_n(r)+ c_{2,n,\omega} Y_n(r)
\end{eqnarray*}
for arbitrary constants $c_{1,n,\omega}$ and $c_{2,n,\omega}$.

It is convenient to switch the limits of integration in the second integral above
and use the fact that $f(t,r,\theta)$ is supported in $|r|<1/2$
to express $\tilde\psi_n$ in the form
\begin{eqnarray}
\tilde \psi_n(\omega,r) &=&\frac{\pi}{2}Y_n(r)\int_0^r \frac{s J_n(s)}{\omega^2}\hat f_n(\omega,\frac{s}{\omega})\,ds +\frac{\pi}{2} J_n(r)\int_r^{\omega/2}\frac{s Y_n(s)}{\omega^2}\hat f_n(\omega,\frac{s}{\omega})\,ds \nonumber\\
&&\hskip 75pt +c_{1,n,\omega} J_n(r)+ c_{2,n,\omega} Y_n(r)\nonumber\\
 &=&\frac{\pi}{2}Y_n(r)\int_0^{r/\omega} J_n(\omega s) s \hat f_n(\omega,s)\,ds + \frac{\pi}{2} J_n(r)\int_{r/\omega}^{1/2}Y_n(\omega s)s\hat f_n(\omega,s)\,ds \nonumber\\
 &&\hskip 75pt +c_{1,n,\omega} J_n(r)+ c_{2,n,\omega} Y_n(r) \label{eqn:stuff37}\\
 &=&Z_{n}(\omega,r) +c_{1,n,\omega} J_n(r)+ c_{2,n,\omega} Y_n(r) ~,\nonumber
\end{eqnarray}
where we make a substitution $s\rightarrow \omega s$ in both integrals in the first line to obtain the integral in
(\ref{eqn:stuff37}).
%\footnote{To see this we note that $J_n$ and $Y_n$ are both solutions to the homogeneous Bessel equation, and then to prove that $Z_n(\omega,r)$ is a particular solution to the inhomogeneous Bessel equation we also use the fact that the Wronskian
%$J_n(r)Y'_n(r)-J'_n(r)Y_n(r)$ equals $\frac{2}{\pi r}$}.
For fixed $\omega>0$ one has the well known asymptotics (see \cite{Watson})
$$
J_0(\omega s) \approx 1 \hbox{ as } s \rightarrow 0 ~~ \hbox{ and } ~~  J_n(\omega s) \approx \frac{1}{n!} \left(\frac{\omega s}{2}\right)^n \hbox{ as } s\rightarrow 0 \hbox{ for any } n\ge 1~,
$$
$$
Y_0(\omega s) \approx \frac{2}{\pi} \log \omega s \hbox{ as } s \rightarrow 0_+ ~~ \hbox{ and } ~~ Y_n(\omega s) \approx -\frac{(n-1)!}{\pi} \left (\frac{\omega s}2\right )^{-n} \hbox{ as } s\rightarrow 0_+ \hbox{ for any } n\ge 1~.
$$
Due to these asymptotic relations the particular solution $Z_n(\omega,r)$ is bounded as $r \rightarrow 0$ for all $n\ge 0$. In order for the function $\tilde \psi_n$ to be bounded as $r \rightarrow 0$ we thus need that the constant $c_{2,n,\omega}$ equals $0$. Finally, to enforce that $\tilde \psi_n$ satisfies the outgoing radiation condition we must have that
$$
c_{1,n,\omega}=-i \frac{\pi}{2}\int_0^{1/2} J_n(\omega s) s \hat f_n(\omega,s)\,ds ~,
$$
because in that case $\tilde \psi_n$ becomes proportional to the H\"ankel function $H_n^{(1)}(r) = J_n(r)+iY_n(r)$ for $r>\omega/2$ (note
that for $r>\omega/2$ the integral multiplied by $J_n(r)$ in (\ref{eqn:stuff37}) equals $0$.)

We can now return to the relation
$\hat \phi_n(\omega,r)=\tilde \psi_n(\omega,\omega r)$ to conclude (for $\omega>0$) that\footnote{The formula (\ref{phift}) only holds for  $\omega>0$; to maintain the identity (\ref{uhatsum}) for all $\omega \in \R$, we define $\hat \phi_n(\omega,r)= \overline{\hat \phi_n(-\omega,r)}$ for $\omega<0$}
\begin{eqnarray}
\label{phift}
&&\hat\phi_n(\omega,r)=\frac{\pi}{2}Y_n(\omega r)\int_0^{r} J_n(\omega s) s \hat f_n(\omega,s)\,ds \\
&& \hskip 75pt +\frac{\pi}{2} J_n(\omega r)\left( \int_{r}^{1/2}Y_n(\omega s)s\hat f_n(\omega,s)\,ds -i\int_0^{1/2} J_n(\omega s) s \hat f_n(\omega,s)\,ds\right)~. \nonumber
\end{eqnarray}
%\marginpar{Hmm, $\hat\phi_0$ can have a log singularity at $\omega=0$ though; for $n\geq 1$ no problem.}
In particular, for $R>1/2$ the integral in (\ref{phift}) with limits $r=R$ to $r=1/2$ vanishes and we have
\begin{eqnarray}
\label{Rphin}
\hat \phi_n(\omega,R) &=& -i \frac{\pi}{2}\left( J_n(\omega R)+i Y_n(\omega R)\right)\int_0^{1/2} J_n(\omega s)s\hat f_n(\omega,s)\,ds  \nonumber\\
&=&-i \frac{\pi}{2}H_n^{(1)}(\omega R) \int_0^{1/2} J_n(\omega s)s\hat f_n(\omega,s)\,ds ~.
\end{eqnarray}
In combination with (\ref{uhatsum}) this now yields
\begin{equation}
\label{seconduhat}
\hat u(\omega,R,\theta)= i \frac{\pi}{2} \sum_{n=0}^N \cos n\theta \,H_n^{(1)}(\omega R) \int_0^{1/2} J_n(\omega s)s\hat f_n(\omega,s)\,ds ~,
\end{equation}
for any $R>1/2$, $\omega>0$.

\begin{remark}

From the above formula and the aforementioned formula for the inversion of the Fourier Transform we get, for $R>1/2$,
\begin{eqnarray*}
u(0,R,\theta)&=&-\Im \sqrt{\frac{\pi}{2}}\sum_{n=0}^N \cos n\theta \,\int_0^\infty H_n^{(1)}(\omega R) \int_0^{1/2} J_n(\omega s)s\hat f_n(\omega,s)\,ds\,d\omega \\
&=&-\sqrt{\frac{\pi}{2}}\sum_{n=0}^N \cos n\theta \int_0^{1/2} s \int_0^\infty J_n(\omega s)\Big( Y_n(\omega R) \Re \hat f_n(\omega,s) \\
&& \hskip 125pt + J_n(\omega R) \Im \hat f_n(\omega,s)\Big) d\omega\, ds~,
\end{eqnarray*}
as well as
\begin{eqnarray*}
\frac{\partial u}{\partial t} (0,R,\theta)&=&\Re \sqrt{\frac{\pi}{2}}\sum_{n=0}^N \cos n\theta \,\int_0^\infty \omega H_n^{(1)}(\omega R) \int_0^{1/2} J_n(\omega s)s\hat f_n(\omega,s)\,ds\,d\omega \\
&=&\sqrt{\frac{\pi}{2}}\sum_{n=0}^N \cos n\theta \int_0^{1/2} s \int_0^\infty \omega J_n(\omega s)\Big( J_n(\omega R) \Re \hat f_n(\omega,s) \\
&& \hskip 125pt - Y_n(\omega R) \Im \hat f_n(\omega,s)\Big) d\omega\, ds~.
\end{eqnarray*}
Now consider the special case $N=0$ and $f_0(t,r)= 1_{\{ 0<t<2L\}}\times 1_{\{r<\frac12\}}$, then
$$
\hat f_0(\omega,r) = \frac{1}{\sqrt{2\pi}}\left( \frac{\sin 2L \omega}{\omega} +i 2 \frac{\sin^2 L\omega }{\omega} \right) ~~ \hbox{ for } r<1/2~.
$$
and $\hat f_0(\omega,r)=0$ for $r>1/2$. Therefore, in this case actually for $R>1/2$,
\begin{eqnarray}
\label{u1}
u(0,R,\theta)&=&-\frac{1}{2} \int_0^{1/2} s \int_0^\infty J_0(\omega s)\Big( Y_0(\omega R) \frac{\sin 2 L \omega}{\omega}  \\
&& \hskip 125pt + 2 J_0(\omega R) \Im \frac{\sin^2L\omega}{\omega}\Big) d\omega\, ds \nonumber \\
&=&-\frac{1}{4} I_1(R) \nonumber
\end{eqnarray}
with
$$I_1(R)=\int_0^\infty \frac{J_1(\omega/2)}{\omega}\Big( Y_0(\omega R) \frac{\sin 2 L \omega}{\omega}  + 2 J_0(\omega R) \frac{\sin^2L\omega}{\omega}\Big) d\omega ~, $$
and, similarly
\begin{equation}
\label{u2}
\frac {\partial u}{\partial t}(0,R,\theta) =\frac{1}{4} I_2(R)
\end{equation}
with
$$I_2(R)= \int_0^\infty J_1(\omega/2)\Big( J_0(\omega R) \frac{\sin 2 L \omega}{\omega} - 2 Y_0(\omega R) \frac{\sin^2L\omega}{\omega}\Big) d\omega ~.$$

The left and right panels in Figure \ref{fig:Vplots} show the graphs of the functions $I_1(R)$ and $I_2(R)$, which are consistent with the fact that $u(0,\cdot) = \frac {\partial u}{\partial t}(0,\cdot)=0$ in our initial value problem for the wave equation, and the fact that the formulas (\ref{u1}-\ref{u2}) hold for $R>1/2$. Note that these integrals do
not represent $u(0,\cdot)$ or $\frac{\partial u}{\partial t}(0,\cdot)$ for $R<1/2$.
\begin{figure}[h]
\centerline{\includegraphics[width=0.5\textwidth]{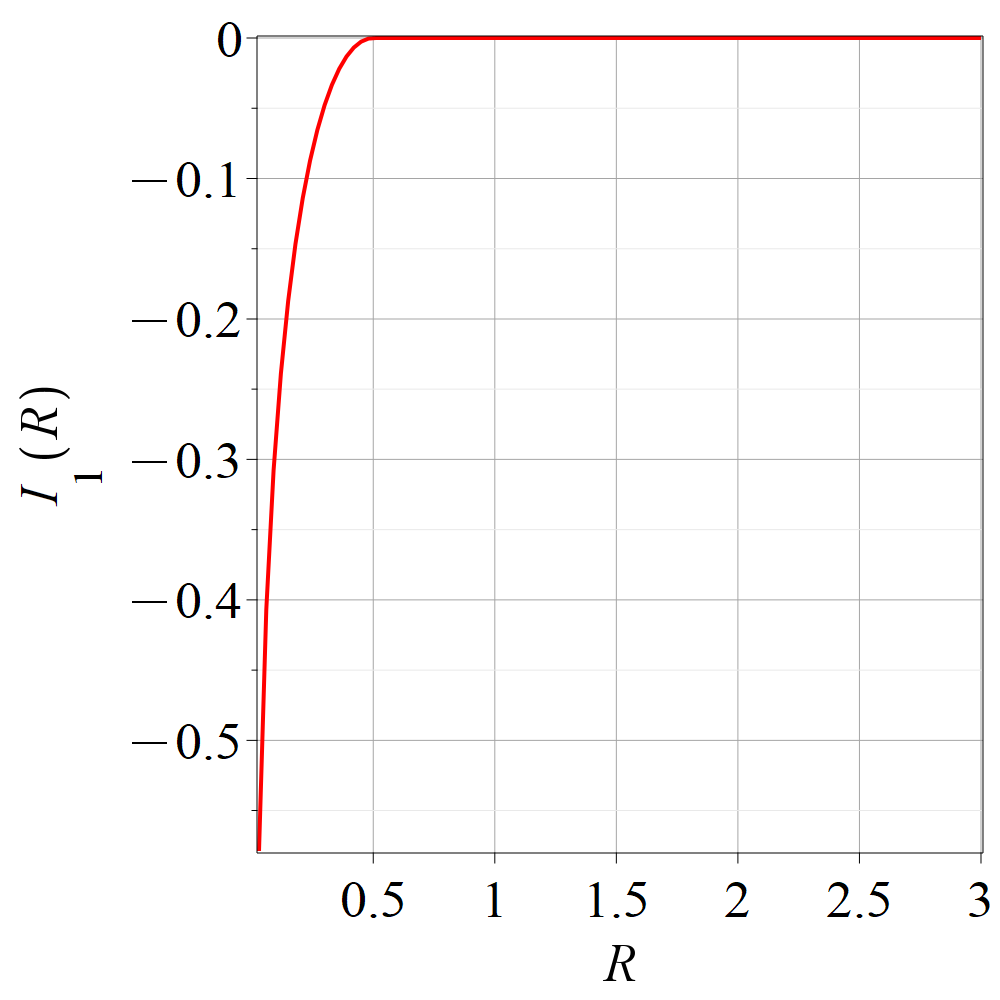}\includegraphics[width=0.5\textwidth]{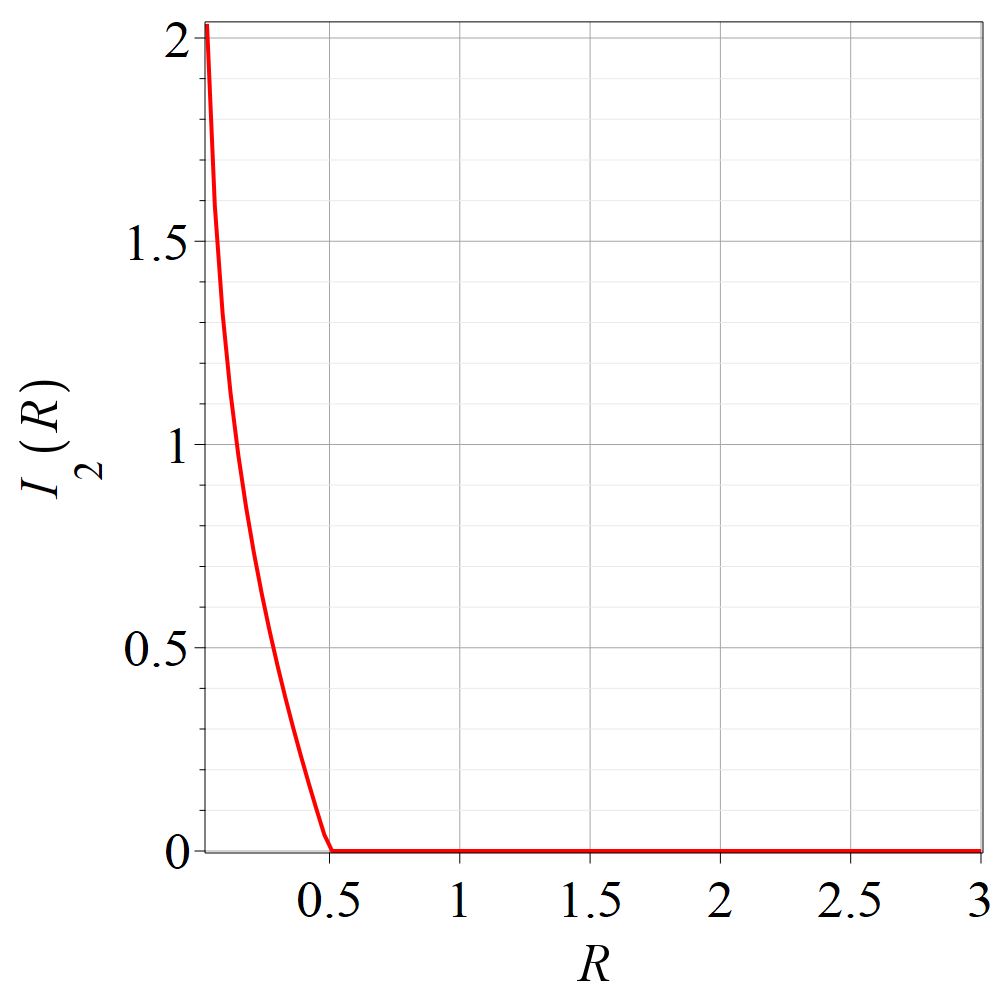}}
\caption{Left panel: graph of $I_1(R)$ for $0\leq R\leq 3$ with $L=3$. Right panel: graph of $I_2(R)$ for $0\leq R\leq 3$ with $L=3$.}
\label{fig:Vplots}
\end{figure}
\end{remark}

%NOTE: The integrals I_1 and I_2 have logarithmic singularities at r = 0.

\section{A PML-modified solution}
\label{sec:pmlsolution}

Let $\alpha(\cdot)$ be a smooth monotone function on $[0,\infty)$, with
\begin{equation}
\label{eqn:alphadef}
\alpha(r)= 0 ~\hbox{ for }~ 0\le r\le 1~~,~~\hbox{ and } \alpha(r)= \alpha_0>0 ~\hbox{ for }~ 2\le r~.
\end{equation}
We define an auxiliary function $V$ associated with the wave equation (\ref{wave-eq}) via its Fourier transform
$\hat{V}$ with respect to $t$ according to
\begin{eqnarray}
\hat V(\omega,r,\theta) & = & \hat u(\omega, r(1+i\alpha(r)/\omega),\theta)\nonumber\\
& = & i \frac{\pi}{2} \sum_{n=0}^N \cos n\theta \,H_n^{(1)}(\omega r\left(1+i\frac{\alpha(r)}{\omega}\right)) \int_0^{1/2} J_n(\omega s)s\hat f_n(\omega,s)\,ds,\label{rless12}
\end{eqnarray}
for $r>1/2$ \footnote{For $1/2<r<1$, or more precisely, when $\alpha(r)=0$, the formula (\ref{rless12}) only applies to $\omega>0$, for $\omega<0$ we define $\hat V(\omega,r,\theta)=\overline{\hat V(-\omega,r,\theta)}$.}, and
\begin{equation}
\label{Vhatsum}
\hat V(\omega,r,\theta)= \hat u(\omega,r,\theta) = -\sum_{n=0}^N \cos n\theta \, \hat \phi_n(\omega,r)
\end{equation}
for $0\le r\le 1/2$, where $\hat \phi_n(\omega,r)$ is given by (\ref{phift}). We have effectively defined
$\hat V$ to equal $\hat u$ as given by (\ref{seconduhat}) but along a deformed path in the complex $r$-plane parameterized by
$r(1+i\alpha(r)/\omega)$ for $r>0$.
Note that for any $\omega>0$ and $r>1/2$ the H\"ankel functions $H^{(1)}_n$
are defined and analytic on this path. Also, given the properties of $\alpha$ and the expression for $\hat\phi_n$ in (\ref{phift})-(\ref{Rphin})
we see from (\ref{rless12})-(\ref{Vhatsum}) that $\hat V$ transitions smoothly across $r=1/2$.
%\marginpar{If $0<r<1$ then $V$ may be logarithmically singular at $\omega=0$.}

For a fixed (large) $R>1/2$, we also define
$$
\hat W(\omega,r,\theta)=  -i\frac{\pi}{2}\sum_{n=0}^N  d_n \cos n\theta J_n (\omega r(1+i\alpha(r)/\omega ))
$$
where the coefficients $d_n$ are given by
\begin{equation}
\label{dnform}
d_n(\omega) = \frac{H_n^{(1)}(\omega R(1+i\alpha(R)/\omega )}{J_n(\omega R(1+i\alpha(R)/\omega )} \int_0^{1/2} J_n(\omega s)s\hat f_n(\omega,s)\,ds~.
\end{equation}

\begin{lemma}
\label{lemma:pmlfunc}
Let $u_{pml}(t,r,\theta)$ be defined as the inverse Fourier transform in time of the function
\begin{equation}
\label{uPML}
\hat u_{pml}=\hat V(\omega,r,\theta)+\hat W(\omega,r,\theta).
\end{equation}
Then for $r<R$ we conclude that $u_{pml}$ is real-valued. We will refer to this as
the ``$R$-truncated PML version'' of the solution to the wave equation (\ref{wave-eq}).
\footnote{In other words the $R$-truncated PML version of the solution to the wave equation (\ref{wave-eq}) is given by $u_{pml}(t,r,\theta) = 2 \Re \left(\frac1{\sqrt{2\pi}} \int_0^\infty e^{-i\omega t} (\hat V(\omega,r,\theta)+\hat W(\omega,r,\theta))\,d\omega \right)$, $r<R$.}
Also, $u_{pml}(t,R,\theta)=0$ for all $t>0$ and $\theta$, and
$$
u_{pml}(0,r,\theta) = \frac{\partial u_{pml}}{\partial t}(0,r,\theta) = 0
$$
for all $r<R$ and $\theta$.
\end{lemma}
{\bf Proof:} From the symmetries stated in Lemma \ref{relcompl} and Lemma \ref{relcompl3} of the Appendix and the symmetries of $\hat \phi_n$ and $d_n$ (see (\ref{cdsym})) we obtain the symmetries
$$
\hat V(-\omega,r,\theta)= \overline{\hat V(\omega,r,\theta)}~~\hbox{ and } ~~\hat W(-\omega,r,\theta)= \overline{\hat W(\omega,r,\theta)}~,
$$
also when $\alpha(r)>0$. We conclude that $u_{pml}$ as given by (\ref{uPML})
is the Fourier Transform of a real valued function
We also note that $\hat u_{pml} (\omega,r,\theta)$ (and thus $u_{pml}(t,r,\theta))$ vanishes at $r=R$, and that according to (\ref{uhatsum})-(\ref{seconduhat}) and (\ref{rless12})-(\ref{Vhatsum})
$\hat u=\hat V$ for $r<1$, since $\alpha(r)$ vanishes there. Finally, contour integration of $\hat W(\omega,r,\theta)$ and $\hat V(\omega,r,\theta)$ (the latter for $r>1/2$) in the complex upper half-plane, with a semi-circular arc expanding to $\infty$, yields that
$$
\int_{-\infty}^{\infty}\hat W(\omega,r,\theta)\,d\omega = 0~,~~\hbox{ for } r>0~, ~~\hbox{and } \int_{-\infty}^{\infty}\hat V(\omega,r,\theta)\,d\omega = 0~,~~\hbox{ for } r>1/2~,
$$
provided the functions $\hat f_n(\cdot,s)$ are analytic and uniformly bounded in the complex upper half-plane (here we also use the asymptotics of the H\"ankel functions given below). This, in combination with the fact that $V=u$ for $r<1$ shows that
$$
V(0,r,\theta)=W(0,r,\theta)=0 ~~ \hbox{ and thus } ~~ u_{pml}(0,r,\theta)=0~ \hbox{ for all }~ r<R ~ \hbox{ and all }~~ \theta~.
$$
Under the additional assumption that the functions $\hat f_n(\omega,s)$ go to zero as $|\omega|^{-1}$ as $|\omega|\rightarrow \infty$ in the complex upper half-plane, one may similarly show that
$$
\frac{\partial}{\partial t}V(0,r,\theta)=\frac{\partial}{\partial t}W(0,r,\theta)=0 ~~ \hbox{ and thus } ~~ \frac{\partial}{\partial t}u_{pml}(0,r,\theta)=0
$$
for all $r<R$ and all $\theta$. This establishes Lemma \ref{lemma:pmlfunc}.\\\ \\
%Checked numerically, did not yet check analytical details.

The difference between the solution to the infinite space wave equation (\ref{wave-eq}) and the $R$-truncated PML version
$u_{pml}$ of Lemma \ref{lemma:pmlfunc} is on the region $0\le r\le 1$ exactly given by
\begin{eqnarray*}
W(t,r,\theta) &=& \sqrt{\frac{2}{\pi}} \Re \int_{0}^\infty e^{-i t \omega}\hat W(\omega,r,\theta)\,d\omega~ \\
&=& \sqrt{\frac{{\pi}}{2}}\sum_{n=0}^N  \cos n\theta \\
&&\hskip 5pt \times \Im \int_{0}^\infty e^{-it\omega} \frac{H_n^{(1)}(\omega R\left(1+i\alpha(R)/\omega\right))}{J_n (\omega R\left(1+i\alpha(R)/\omega\right))} J_n (\omega r) \int_0^{1/2} J_n(\omega s)s\hat f_n(\omega,s)\,ds \,d\omega~.
\end{eqnarray*}
If we furthermore assume that $R>2$ then this simplifies to
\begin{eqnarray}
W(t,r,\theta) &=&\sqrt{\frac{\pi}{2}}\sum_{n=0}^N  \cos n\theta \nonumber\\
 && \hskip 5pt \times \Im \int_{0}^\infty e^{-it\omega} \frac{H_n^{(1)}(\left(\omega+i\alpha_0\right)R)}{J_n (\left(\omega+i\alpha_0\right)R)} J_n (\omega r) \int_0^{1/2} J_n(\omega s)s\hat f_n(\omega,s)\,ds \,d\omega~.\label{eqn:Wdef2}
\end{eqnarray}
We have the following asymptotic results for the H\"ankel functions of the first and second kind
$$
H_n^{(1)}(z)=\left(\frac{2}{\pi z}\right)^{\frac12}e^{i(z-\frac{n}2-\frac14 \pi)}\left( 1+O(\frac{1}{|z|})\right)~,~-\pi<\mbox{arg}(z)<\pi,
$$
and
$$
H_n^{(2)}(z)=\left(\frac{2}{\pi z}\right)^{\frac12}e^{-i(z-\frac{n}2-\frac14 \pi)}\left( 1+O(\frac{1}{|z|})\right)~,~ -\pi<\mbox{arg}(z)<\pi.
$$
These asymptotic formulae are valid uniformly for $\mbox{arg}(z)\in (-\pi+\delta,\pi-\delta)$ for any $\delta>0$ (see \cite{Watson} section 7.2) In particular, if we consider $z$ of the form $z=(\omega+i\alpha_0)R$ with $\alpha_0>0$ fixed, then
$$
H_n^{(1)}((\omega +i\alpha_0)R)=\left(\frac{2}{\pi (\omega+i \alpha_0)}\right)^{\frac12}R^{-1/2}e^{-\alpha_0\,R +i \omega R-i(\frac{n}2+\frac14 \pi)}\left(1+O(\frac{1}{\alpha_0 R})\right)~,
$$
uniformly in $\omega\ge 0$, $R>1$, and
$$
H_n^{(2)}((\omega +i\alpha_0)R)=\left(\frac{2}{\pi (\omega+i \alpha_0)}\right)^{\frac12}R^{-1/2}e^{\alpha_0\,R -i \omega R+i(\frac{n}2+\frac14 \pi)}\left(1+O(\frac{1}{\alpha_0 R})\right)~,
$$
uniformly in $\omega\ge 0$, $R>1$.
In combination with the fact that
$$
J_n(z)=\frac12\left(H_n^{(1)}(z)+H_n^{(2)}(z) \right)~,
$$
we thus arrive at
\begin{lemma}
\label{asympest}
The ratio $\frac{H_n^{(1)}(\left(\omega+i\alpha_0\right)R)}{J_n (\left(\omega+i\alpha_0\right)R)}$ equals $2e^{-2\alpha_0\,R +i 2\omega R-i(n+\frac{\pi}2)}\left(1+O(1/(\alpha_0R)\right)$, where the term $O(1/(\alpha_0R))$ converges to zero as  $R \rightarrow \infty$ uniformly with respect to $\omega\ge 0$ for any fixed real number $\alpha_0>0$ and any fixed integer $n\ge 0$. As a consequence, given any fixed $\alpha_0>0$ and $n\ge 0$ there exist positive constants $R_0$, $c$ and $C$ (independent of $\omega$ and $R$) such that
$$
c e^{-2\alpha_0\,R}\le \left|\frac{H_n^{(1)}(\left(\omega+i\alpha_0\right)R)}{J_n (\left(\omega+i\alpha_0\right)R)}\right|\le Ce^{-2\alpha_0\,R}~,
$$
for all $\omega\ge 0$ and $R>R_0$.
\end{lemma}
\vskip 5pt
\noindent

%\begin{proposition}
%\label{fracest}
%Given any fixed $\alpha_0>0$ and any fixed non-negative integer $n$ , there exist positive constants $R_0$, $c$ and %$C$ (independent of $\omega$ and $R$) such that
%$$
%c e^{-2\alpha_0\,R}\le \left|\frac{H_n^{(1)}(\left(\omega+i\alpha_0\right)R)}{J_n (\left(\omega+i\alpha_0\right)R)}\right|\le %Ce^{-2\alpha_0\,R}~,
%$$
%for all $\omega \in \R$ and $R>R_0$.
%\end{proposition}

\begin{remark}
Due to the continuity and non-vanishing of the expression $\left|\frac{H_n^{(1)}(\left(\omega+i\alpha_0\right)R)}{J_n (\left(\omega+i\alpha_0\right)R)}\right|$ it is very easy to see that we may take $R_0$ to be any fixed positive number, for example $R_0=2$.
\end{remark}

\section{An example}
\label{sec:example}

For a fixed $L>0$ let us take $f_n(t,r) = a_n 1_{\{ 0<t<2L\}}\times 1_{\{r<\frac12\}}$. In that case
\begin{eqnarray}
\label{eqn:fncomp}
\hat f_n(\omega,r) &=&a_n \sqrt{\frac{2}{\pi}}\frac{e^{iL\omega}}{\omega}\sin L \omega \\
&=& a_n \frac{1}{\sqrt{2\pi}} \frac{\sin 2L \omega}{\omega} +i a_n \sqrt{\frac{2}{\pi}} \frac{\sin^2 L\omega }{\omega}\hbox{ for } r<1/2~. \nonumber
\end{eqnarray}
and $\hat f_n(\omega,r)=0$ for $r>1/2$. If we also assume that $r<1$ and  $R>2$ then $W(t,r,\theta)$ defined
by (\ref{eqn:Wdef2}) takes the form
\begin{eqnarray*}
W(t,r,\theta)&=& \sum_{n=0}^N  a_n \cos n\theta \\
&&\hskip 5pt \times \Im \int_{0}^\infty \frac{e^{i(L-t)\omega}}{\omega}\sin L\omega \frac{H_n^{(1)}(\omega R\left(1+i\frac{\alpha_0}{\omega}\right))}{J_n (\omega R\left(1+i\frac{\alpha_0}{\omega}\right))} J_n (\omega r) \int_0^{1/2} J_n(\omega s)s\,ds \,d\omega
\end{eqnarray*}
if we note that $\hat f_n(\omega,r)$ in the first line of (\ref{eqn:fncomp}) is independent of $r$, so that
$\hat f_n(\omega,s)$ may be pulled out of the $ds$ integral above.

We now proceed to estimate this $W(t,r,\theta)$, for $r<1$,  using the estimate in Proposition \ref{asympest}.
\begin{eqnarray}
\label{West}
|W(t,r,\theta)|&\le&  Ce^{-2\alpha_0\,R} \sum_{n=0}^N |a_n| \int_{0}^{\infty} \frac{|J_n(\omega r)|}{\omega}\left|\sin L \omega\right|\left|\int_0^{1/2} J_n(\omega s) s\,ds\right|\,d\omega  \nonumber\\
&\le&Ce^{-2\alpha_0\,R} \sum_{n=0}^N |a_n| \int_{0}^{\infty} \frac{|J_n(\omega r)|}{\omega}\min\{1, L\omega \}\left|\int_0^{1/2} J_n(\omega s) s\,ds\right|\,d\omega \nonumber\\
&=& Ce^{-2\alpha_0\,R} \sum_{n=0}^N |a_n| L\int_{0}^{1/L} |J_n(\omega r)|\left|\int_0^{1/2} J_n(\omega s) s\,ds\right|\,d\omega\\
&& \hskip 10pt + Ce^{-2\alpha_0\,R} \sum_{n=0}^N |a_n| \int_{1/L}^{\infty}\frac{|J_n(\omega r)|}{\omega}\left|\int_0^{1/2} J_n(\omega s) s\,ds\right|\,d\omega~. \nonumber
\end{eqnarray}
A simple calculation gives that
$$
\left| \int_0^{1/2}J_n(\omega s) s\,ds\right| =\omega^{-2}\left| \int_0^{\omega/2}J_n(s)s\,ds \right|\\~,
$$
for $\omega>0$.
As $J_n(x)$ is asymptotic to $\frac{x^n}{2^n n!}$ as $x \rightarrow 0$, and $J_n(x)$ is  bounded by $C x^{-1/2}$ for $x>1$ it follows that
$$
\left| \int_0^{\omega/2}J_n(s)s\,ds \right| \le C \min \{ \omega^{n+2}, \omega^{3/2}  \}
$$
so that, altogether
$$
\left| \int_0^{1/2}J_n(\omega s) s\,ds\right| \le C \min \{ \omega^{n}, \omega^{-1/2}  \}\le C \min \{ 1, \omega^{-1/2}  \}.
$$
From these estimates it follows immediately that
\begin{equation}
\label{aux1}
L\int_{0}^{1/L} |J_n(\omega r)|\left|\int_0^{1/2} J_n(\omega s) s\,ds\right|\,d\omega \le C ~, ~~~ \hbox{ for } 0\le n\le N~,
\end{equation}
and
\begin{eqnarray}
\label{aux2}
&&\int_{1/L}^{\infty}\frac{|J_n(\omega r)|}{\omega}\left|\int_0^{1/2} J_n(\omega s) s\,ds\right|\,d\omega \nonumber \\
&& \hskip 30pt = \int_{1/L}^{\max\{1,1/L\}}\frac{|J_n(\omega r)|}{\omega}\left|\int_0^{1/2} J_n(\omega s) s\,ds\right|\,d\omega \nonumber \\
&&\hskip 50pt + \int_{\max\{1,1/L\}}^{\infty}\frac{|J_n(\omega r)|}{\omega}\left|\int_0^{1/2} J_n(\omega s) s\,ds\right|\,d\omega \nonumber\\
&&\hskip 30pt \le C \int_0^1 \omega^{n-1}\,d\omega +C\int_1^\infty \omega^{-3/2}\,d\omega \le C~, ~~~ \hbox{ for } 1\le n \le N~,
\end{eqnarray}
whereas for $n=0$
\begin{eqnarray}
\label{aux3}
&&\int_{1/L}^{\infty}\frac{|J_0(\omega r)|}{\omega}\left|\int_0^{1/2} J_0(\omega s) s\,ds\right|\,d\omega \nonumber \\
&& \hskip 30pt = \int_{1/L}^{\max\{1,1/L\}}\omega^{-1}\,d\omega \nonumber + \int_{\max\{1,1/L\}}^{\infty} \omega^{-3/2}\nonumber \\
&&\hskip 30pt \le C (|\log L |+1)~.
\end{eqnarray}

%\begin{remark}
%We note that for $n=0$ this estimate may very simply be improved to
%$$
%\left| \int_0^{1/2}J_0(\omega s) s\,ds\right| \le C \min \{ 1, |\omega|^{-3/2}  \}
%$$
%by using that
%$$\frac{d}{ds}(s J_1(s))= sJ_0(s).$$
%\end{remark}

Inserting (\ref{aux1})-(\ref{aux3}) into the estimate (\ref{West})we now obtain
$$
|W(t,r,\theta)| \le Ce^{-2\alpha_0\,R} \left( \sum_{n=1}^N |a_n| + (|\log L |+1)|a_0|\right)
$$
for $r<1$, $t>0$ and $\theta \in [0,2\pi]$. The constant $C$ depends on $\alpha_0$ and $N$, but  is independent of $r,t,\theta$, $R$ and $a_n$, { \bf as well as } $L$. In summary we have thus proven
\begin{proposition}
\label{prop:mainresult}
Let $f(t,r,\theta)= f_L(t,r)\sum_{n=0}^N a_n \cos(n\theta)$ with $f_L(t,r)= 1_{\{ 0<t<2L\}}\times 1_{\{r<\frac12\}}$, and
let $u$ be the solution to the inhomogeneous full space wave equation (\ref{wave-eq}) with this right hand side. If
$u_{pml}$ denotes the $R$-truncated PML solution, whose Fourier transform is given by (\ref{uPML}) for some fixed $\alpha_0>0$, then
$$
\max_{\begin{array}{cc} &0\le t,~0\le r< 1,\\ &0\le  \theta<2\pi \end{array}}|u(t,r,\theta)-u_{pml}(t,r.\theta)| \le Ce^{-2\alpha_0\,R}\left(\sum_{n=1}^N |a_n|+ (|\log L |+1)|a_0| \right).
$$
The constant $C$ is independent of $R$, $\{a_n\}$ and $L$.
\end{proposition}

\begin{remark}
\label{rem:mainresult2}
A similar estimation would hold for any function $g(t,r)$ (in place of $f_L(t,r)$) whose Fourier transform is analytic in the complex upper half-plane and goes to zero as $|\omega|^{-1}$ when $\omega \rightarrow \infty$ in the same half-plane.
\end{remark}

\section{Representation of the PML equations as a system in time-domain}
\label{sec:pmlpdes}

For any fixed $n$ and $\omega$ let $K_n(t,\theta,\tilde r)$ denote the function
$$
\tilde K_n(t, \omega,\theta,\tilde r)= e^{-i t\omega }\cos n\theta H_n^{(1)}(\tilde r).
$$
It is well known that $\tilde K_n$ satisfies the equation
\begin{equation}
\label{Keq}
\left[ \frac{1}{\omega^2}\frac{\partial^2}{\partial t^2} -\frac{1}{\tilde r}\frac{\partial}{\partial \tilde r} \left( \tilde r \frac{\partial}{\partial \tilde r} \right) -\frac{1}{\tilde r^2} \frac{\partial ^2}{\partial \theta^2}\right] \tilde K_n =0~.
\end{equation}
Let use define
$$
\tilde r=\omega r(1+i\alpha(r)/\omega)\quad\textrm{and}\quad
\beta(r)=\frac{d}{dr}(r\alpha(r)).
$$
Noting that $\frac{\partial}{\partial \tilde r}= \frac{1}{\omega + i \beta(r)}\frac{\partial}{\partial r}$ we get after insertion into (\ref{Keq}) and multiplication by $\tilde r^2$
$$
\left[ \left (r+i \frac{r \alpha(r)}{\omega}\right )^2\frac{\partial^2}{\partial t^2} -\frac{r+i \frac{r \alpha(r)}{\omega}}{1+i \frac{\beta(r)}{\omega}}\frac{\partial}{\partial  r} \left( \frac{r+i \frac{r \alpha(r)}{\omega}}{1+i \frac{\beta(r)}{\omega}}\frac{\partial}{\partial r} \right) -\frac{\partial ^2}{\partial \theta^2}\right] K_n =0~,
$$
where $K_n(t,\theta,r)$ denotes the function
$$
K_n(t, \omega,\theta, r)= e^{-i t\omega }\cos n\theta H_n^{(1)}(\omega r(1+i\alpha(r)/\omega)).
$$

We introduce auxiliary functions $\Phi_n$ and $\Psi_n$ by
$$
\Phi_n(t,\omega,\theta,r)=  \frac{r+i \frac{r \alpha(r)}{\omega}}{1+i \frac{\beta(r)}{\omega}}\frac{\partial}{\partial r} K_n(t,\omega,\theta,r)
$$
and
$$
\Psi_n(t,\omega,\theta,r)=\frac{r+i \frac{r \alpha(r)}{\omega}}{1+i \frac{\beta(r)}{\omega}}\frac{\partial}{\partial r} \Phi_n(t,\omega,\theta,r)
$$
It follows immediately that these functions satisfy the equations
\begin{equation}
\label{Phieq}
\frac{\partial}{\partial t} \Phi_n +\beta(r) \Phi_n = r\frac{\partial^2}{\partial t \partial r}K_n + r\alpha(r)\frac{\partial}{\partial r}K_n
\end{equation}
and
\begin{equation}
\label{Psieq}
\frac{\partial}{\partial t} \Psi_n +\beta(r) \Psi_n = r\frac{\partial^2}{\partial t \partial r}\Phi_n + r\alpha(r)\frac{\partial}{\partial r}\Phi_n
\end{equation}
It's easy to check that %\marginpar{The $2r^2\alpha(r)\partial K_n/\partial t$ is now right, but was wrong (it was $2r\alpha(r)\partial K_n/\partial t$)}
$$
\left (r+i \frac{r \alpha(r)}{\omega}\right )^2\frac{\partial^2}{\partial t^2} K_n = r^2 \frac{\partial^2}{\partial t^2}K_n + 2r^2\alpha(r) \frac{\partial}{\partial t}K_n+(r\alpha(r))^2K_n,
$$
and so it follows that
\begin{equation}
\label{2ndKeq}
r^2 \frac{\partial^2}{\partial t^2}K_n +2r^2\alpha(r) \frac{\partial}{\partial t}K_n+(r\alpha(r))^2K_n -\frac{\partial^2}{\partial \theta^2}K_n =\Psi_n .
\end{equation}
The same set of equations are satisfied by the functions $K_n^0$, $\Phi_n^0$ and $\Psi_n^0$ that arise when  $H_n^{(1)}$ is replaced by $J_n$ in the above definitions.

We note that for $r>1/2$
$$
e^{-i\omega t}\hat u_{pml} =i\frac{\pi}{2}\sum_{n=0}^N \left(c_n K_n(t,\omega,r,\theta) - d_n K_n^0(t,\omega,r,\theta)\right)~,
$$
where $d_n(\omega)$ is given by (\ref{dnform}) and $c_n(\omega)$ is given by
$$
c_n(\omega)=\int_0^{1/2}J_n(\omega s) s \hat f_n(\omega,s)\,ds~.
$$
For $r>1/2$ we introduce functions $v_{pml}$ and $w_{pml}$ by means of the formulas
\begin{equation}
\label{def1}
e^{-i\omega t}\hat v_{pml}=i\frac{\pi}{2}\sum_{n=0}^N\frac{r+i\frac{r\alpha(r)}{\omega} }{1+i\frac{\beta(r)}{\omega}}\frac{\partial}{\partial r}\left(c_n K_n(t,\omega,r,\theta) - d_n K_n^0(t,\omega,r,\theta)\right)
\end{equation}
and
\begin{eqnarray}
\label{def2}
e^{-i\omega t}\hat w_{pml}&=&i\frac{\pi}{2}\sum_{n=0}^N \frac{r+i\frac{r\alpha(r)}{\omega} }{1+i\frac{\beta(r)}{\omega}}\frac{\partial}{\partial r}\left(\frac{r+i\frac{r\alpha(r)}{\omega} }{1+i\frac{\beta(r)}{\omega}}\frac{\partial}{\partial r} \right) \left(c_n K_n(t,\omega,r,\theta) - d_n K_n^0(t,\omega,r,\theta)\right) \nonumber\\
&=&i\frac{\pi}{2}\sum_{n=0}^N \left( n^2 - (\omega r+i r \alpha(r))^2  \right)\left(c_n K_n(t,\omega,r,\theta) - d_n K_n^0(t,\omega,r,\theta)\right)
\end{eqnarray}
Based on the symmetries in Lemma \ref{relcompl} and Lemma \ref{relcompl3} of the Appendix we conclude that
\begin{equation}
\label{cdsym}
\overline{c_n(\omega)}= (-1)^n c_n(-\omega) ~~ \hbox{ and } ~~ \overline{d_n(\omega)}= (-1)^{n+1} d_n(-\omega) ~~ \hbox{ for } \omega \in \R~.
\end{equation}
Using the identity
$$
\frac{r+i\frac{r\alpha(r)}{\omega} }{1+i\frac{\beta(r)}{\omega}}\frac{\partial}{\partial r}H_n^{(1)}(\omega r(1+i\frac{\alpha(r)}{\omega}))= \omega r(1+i\frac{\alpha(r)}{\omega}) \frac{d}{dz}H_n^{(1)}(\omega r(1+i\frac{\alpha(r)}{\omega}))
$$
and a similar identity for $J_n$, in place of $H_n^{(1)}$, in combination with (\ref{cdsym}) and the symmetries in Lemma \ref{relcompl} and Lemma \ref{relcompl2} we now arrive at the fact that
$$
\overline{\hat v_{pml}(\omega)}=\hat v_{pml}(-\omega) ~~ \hbox{ and } ~~ \overline{\hat w_{pml}(\omega)}=\hat w_{pml}(-\omega)~;
$$
in other words: the functions $v_{pml}$ and $w_{pml}$ are also real. After integration with respect to $\omega$ we end up with the following equations for $u_{pml}$ and the associated auxiliary functions $v_{pml}$ and $w_{pml}$ for $r>1/2$:
$$
r^2 \frac{\partial^2}{\partial t^2}u_{pml} +2r^2\alpha(r) \frac{\partial}{\partial t}u_{pml}+(r\alpha(r))^2u_{pml} -\frac{\partial^2}{\partial \theta^2}u_{pml} = w_{pml} ,
$$
\begin{equation}
\label{vpmleq}
\frac{\partial}{\partial t} v_{pml} +\beta(r) v_{pml} = r\frac{\partial^2}{\partial t \partial r}u_{pml}+r\alpha(r)\frac{\partial}{\partial r}u_{pml},
\end{equation}
and
\begin{equation}
\label{wpmleq}
\frac{\partial}{\partial t} w_{pml} +\beta(r) w_{pml} = r\frac{\partial^2}{\partial t \partial r}v_{pml} +r\alpha(r)\frac{\partial}{\partial r} v_{pml}.
\end{equation}
For $r<1/2$ a source term $f(t,r)$ appears on the right hand side of the original wave equation. Since  $\alpha(r)$ and $\beta(r)$ both vanish for $r<1$, it is quite easy to see that $u_{pml}$, $v_{pml}=r\frac{\partial}{\partial r}u_{pml}$ and $w_{pml}=r\frac{\partial}{\partial r}r\frac{\partial}{\partial r}u_{pml}$ satisfy the equation
\begin{equation}
\label{upmleq}
r^2 \frac{\partial^2}{\partial t^2}u_{pml} +2r^2\alpha(r) \frac{\partial}{\partial t}u_{pml}+(r\alpha(r))^2u_{pml} -\frac{\partial^2}{\partial \theta^2}u_{pml} = w_{pml} +r^2 f
\end{equation}
together with the two equations (\ref{vpmleq}) and (\ref{wpmleq}) for $r<1/2$ (actually for $r<1$). We also note that for $1/2<r<1$ $v_{pml}=r\frac{\partial}{\partial r}u_{pml}$  and $w_{pml}=r\frac{\partial}{\partial r}r\frac{\partial}{\partial r}u_{pml}$ is consistent with the definitions (\ref{def1}) and (\ref{def2}).  In other words: given a source $f$ with support inside $\{~r<1/2~\}$ the $R$-truncated PML solution is the first component of a solution to the system (\ref{vpmleq})-(\ref{upmleq}), $r<R$. In the following section we determine the appropriate initial (and boundary) conditions for $u_{pml}$, $v_{pml}$ and $w_{pml}$.
%\marginpar{Should we say anything about well-posedness?}

\section{Initial and boundary conditions for the PML system (\ref{vpmleq})-(\ref{upmleq}) in time-domain}
\label{sec:initbound}

We earlier saw that $u_{pml}$ satisfies the initial conditions
$$
u_{pml}(0,r,\theta)= \frac{\partial}{\partial t}u_{pml}(0,r,\theta)=0~~ \hbox{ for all } ~~r<R ~~\hbox{ and all } ~~ \theta~.
$$
It also satisfies the boundary condition $u_{pml}(t,R,\theta)=0$ for all $t>0$ and $\theta$. Furthermore, $u_{pml}(t,r,\theta)$ has a finite limit as $r \rightarrow 0$, and this limit is independent of $\theta$.

For $r<1$ the initial values $v_{pml}(0,r,\theta)$ and $w_{pml}(0,r,\theta)$ equal zero, due to the initial conditions on $u_{pml}$ and the fact that  $v_{pml}=r\frac{\partial}{\partial r}u_{pml}$  and $w_{pml}=r\frac{\partial}{\partial r}r\frac{\partial}{\partial r}u_{pml}$. For $r>1$ (actually $r>1/2$) a contour integration argument, just as that outlined for $u_{pml}$, gives that
$$
\int_{-\infty}^\infty \hat v_{pml}(\omega,r,\theta)\, d\omega = \int_{-\infty}^\infty \hat w_{pml}(\omega,r,\theta)\, d\omega =0~,
$$
and thus the initial data $v_{pml}(0,r,\theta)$ and $w_{pml}(0,r,\theta)$ vanish for $r>1$ as well.
We note that at $r=0$, $v_{pml}$ and $w_{pml}$ automatically become $0$ since they are constant in $t$ (and have an initial value $0$). All three functions $u_{pml}$, $v_{pml}$ and $w_{pml}$ are $2\pi$-periodic in $\theta$.

%% file: numerics.tex
\section{Numerical Examples}
\label{sec:numerics}

In this section we give a few numerical examples to explore the bound in Proposition \ref{prop:mainresult}, or
more generally, Remark \ref{rem:mainresult2} that follows that proposition. Our goal is not a formal treatment of numerical methods for solving the PML system (\ref{vpmleq})-(\ref{upmleq}), but rather a few convincing numerical computations to illustrate the efficacy of the PML approach and the accuracy of the time domain estimates in Section \ref{sec:example}.

In what follows we work on the disk $D_R=\{(r,\theta);r<R\}$ in the plane with values of $R$ between $R=1$ and $R=3$, and numerically compute the solution to the PML system (\ref{vpmleq})-(\ref{upmleq}) using zero initial data as detailed in
Section \ref{sec:initbound}. We use a zero Dirichlet boundary condition for each of $u_{pml},v_{pml}$, and $w_{pml}$ at $r=R$
(although we should note that the authors in \cite{Monk} remark that other boundary conditions may give better results).
%\marginpar{remark that others (Monk) find other boundary conditions useful.}
The PML function $\alpha(r)$ will be taken as
\[
\alpha(r) = \begin{cases}
0, & 0\leq r \leq 1\\
\alpha_0((3(r-1)^2-2(r-1)^3), & 1<r\leq 2\\
\alpha_0, & r > 2
\end{cases}
\]
which is in accord with (\ref{eqn:alphadef}). We use $\alpha_0=2$.
%Our goal is to compare the
%solution to the PML system (\ref{vpmleq})-(\ref{upmleq}) with the solution to the full-space wave equation (\ref{wave-eq})
%on the region $r<1$.

Computations were performed using a finite element method (FEM) implemented in the
FEniCSx software; see \cite{fenicsx}. A very brief outline of out method is this: we split (\ref{upmleq}) into a pair of first-order equations
in time by introducing an auxiliary function $q=\partial u_{pml}/\partial t$ and so obtain a coupled system of four PDEs, all
first-order in time. We construct a weak form of this system
and then discretize with respect to time. The resulting equations are marched out in time to estimate the solution at times
$t_n=n\Delta t$ where $\Delta t$ is a fixed time step. This is done by
approximating time derivatives with a backward Euler approximation, e.g., $\frac{\partial q}{\partial t}(t_n) \approx (q^{n+1}-q^n)/\Delta t$, where $q^n$ indicates an approximation to $q(t_n)$. The computations are carried out in FEniCSx
using a basis of approximately $15,000$ quadratic Lagrange elements.

\subsection{Analytical solutions to the wave equation on $\mathbb{R}^2$}

Our goal in the numerical experiments that follow is to compare the numerical solution to the PML system
to a numerical evaluation of the exact solution of the inhomogeneous wave equation on $\mathbb{R}^2$ with the same forcing function $f$. In the particular cases of interest that follow we take
\begin{equation}
\label{eqn:gausssource}
f(t,r,\theta) = \frac{50}{\pi}e^{-50r^2} 1_{\{0<t<2L\}} 1_{\{r<1/2\}}
\end{equation}
for some value of $L$, or translations of this function in the plane; note that $f$ here is independent of
$\theta$.
The truncation at $r=1/2$ has little effect, since $f$ decays to approximately $3.7\times 10^{-6}$ at $r=1/2$.
Moreover, since our evaluation of the solution will be at times $t<2L$, the value of $L$ won't matter in these computations
and we may as well assume $L=\infty$.

As a proxy for the solution to the inhomogeneous wave equation (\ref{wave-eq}) on $\R^2$ we work on a disk
of radius $R_2=4$, which is sufficient for times $t\leq 5$ since the wave speed here is $1$ and our interest
is the solution on the unit disk (reflections caused by boundary conditions imposed at $r=4$ do not have time
to return the unit disk). A standard separation of variables for a source term $f=f(t,r)$
shows that the solution to (\ref{wave-eq}) can be expressed as
\begin{equation}
\label{eqn:wexpan}
u(t,r,\theta) = \sum_{n=1}^{\infty} c_n(t) J_0(r\mu_n/R_2)
\end{equation}
where $\mu_n$ is the $n$th positive root of the Bessel function $J_0$ and $c_n(t)$ satisfies
\begin{equation}
\label{eqn:ckode}
c''_n(t) + (\mu_n/R_2)^2 c_n(t) = a_n(t)
\end{equation}
with $c_n(0)=c'_n(0)=0$ for each $n$ and
\begin{equation}
\label{eqn:fntdef}
a_n(t) = \frac{\displaystyle \int_0^{R_2} rf(t,r)J_0(\mu_n r/R_2)\,dr}{\displaystyle  \int_0^{R_2} rJ^2_0(\mu_n r/R_2)\,dr}
\end{equation}
are the coefficients of $f(t,\cdot)$ with respect to the orthogonal basis $J_0(\mu_n r/R_2)$ on the interval $0\leq r\leq R_2$.
With the choice (\ref{eqn:gausssource}) (and take $L=\infty$ so $f$ is independent of $t$) we find $a_n(t)=a_n$ is constant and then
\begin{equation}
\label{eqn:hndef}
c_n(t) = \frac{a_n R_2^2(1-\cos(t\mu_n/R_2))}{\mu_n^2}.
\end{equation}
With $c_n(t)$ as in (\ref{eqn:hndef}), (\ref{eqn:wexpan}) provides the solution $u=u(t,r)$.

\subsection{Experiment 1}

\begin{figure}[h!]
\centerline{\includegraphics[width=0.7\textwidth]{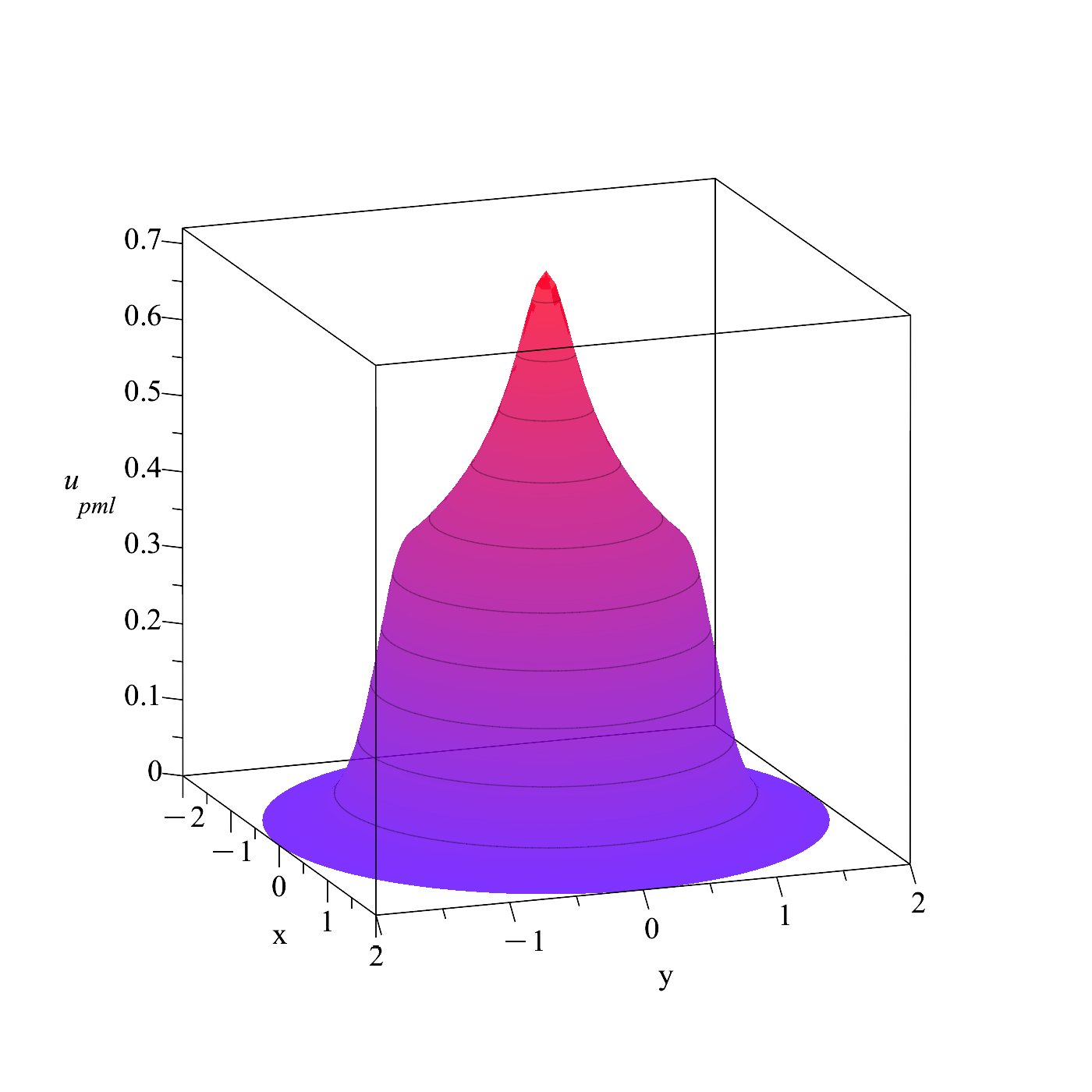}}
\caption{Solution to PML system at time $t=5$ on region $0\leq r\leq 2$}
\label{fig:pmlsol}
\end{figure}

\begin{figure}[h]
\centerline{\includegraphics[width=0.5\textwidth]{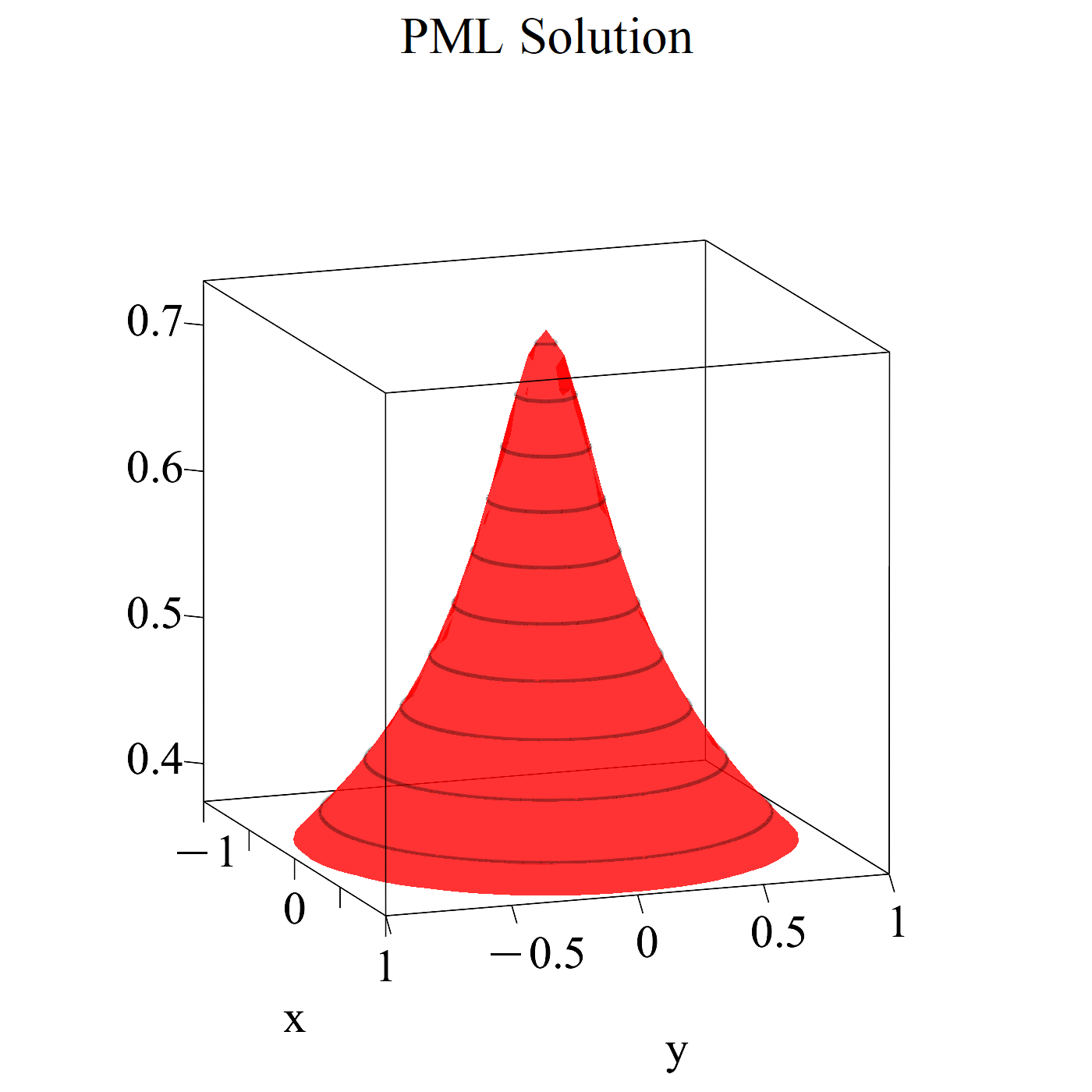}\includegraphics[width=0.5\textwidth]{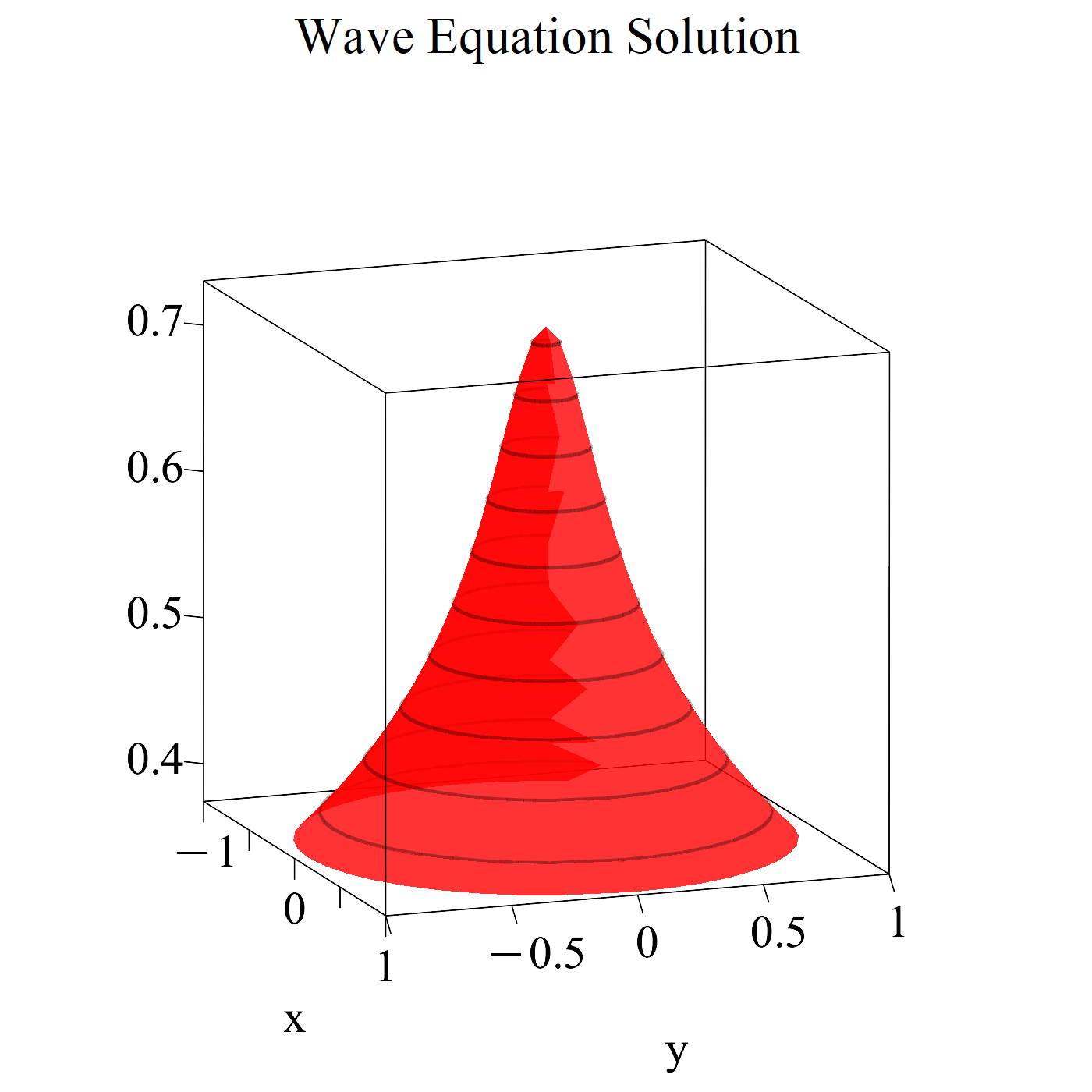}}
\caption{Left panel: Solution to PML system at time $t=5$ on region $0\leq r\leq 1$. Right panel:
Solution to inhomogeneous wave equation on $\mathbb{R}^2$ at time $t=5$ on region $0\leq r\leq 1$.}
\label{fig:pmlsol2}
\end{figure}

We begin by presenting in Figure \ref{fig:pmlsol} a graph of the solution to the PML system (\ref{vpmleq})-(\ref{upmleq}) computed on the region $r<2$. The solution is at time $t=5$ (so without PML any reflections from the boundary at $r=2$ would have plenty of time to return).
Note that the solution decays rapidly as $r$ increases from $r=1$ to $r=2$. The goal of the PML is that the solution
to this system on the region $r<1$ agrees closely with the solution to the corresponding inhomogeneous wave equation
(\ref{wave-eq}) on $r<1$ and this is illustrated in Figure \ref{fig:pmlsol2}. The PML solution and wave equation solutions at $t=5$ (the latter computed using (\ref{eqn:wexpan})-(\ref{eqn:hndef})) are graphed in the left and right panels of Figure \ref{fig:pmlsol2}, respectively.

In this example we find
\[
\max_{0\leq r\leq 1,0\leq \theta\leq 2\pi} |u(5,r,\theta)-u_{pml}(5,r,\theta)| \approx 3.16\times 10^{-3}
\]
while the solution $u$ itself has a supremum norm of approximately $0.72$ at time $t=5$, a relative error
of approximately $0.004$.

\subsection{Experiment 2}

\begin{figure}[h]
\centerline{\includegraphics[width=0.5\textwidth]{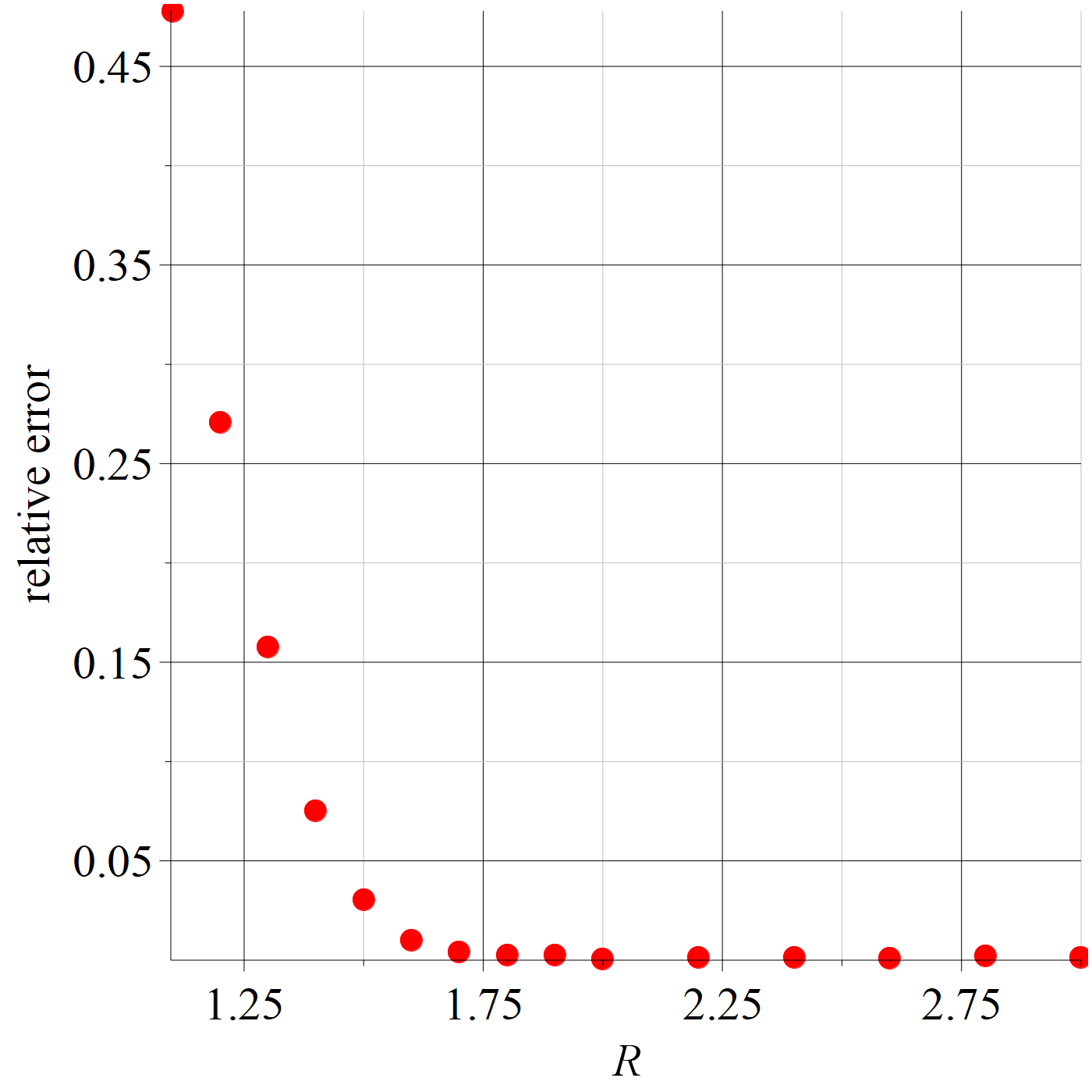}\includegraphics[width=0.5\textwidth]{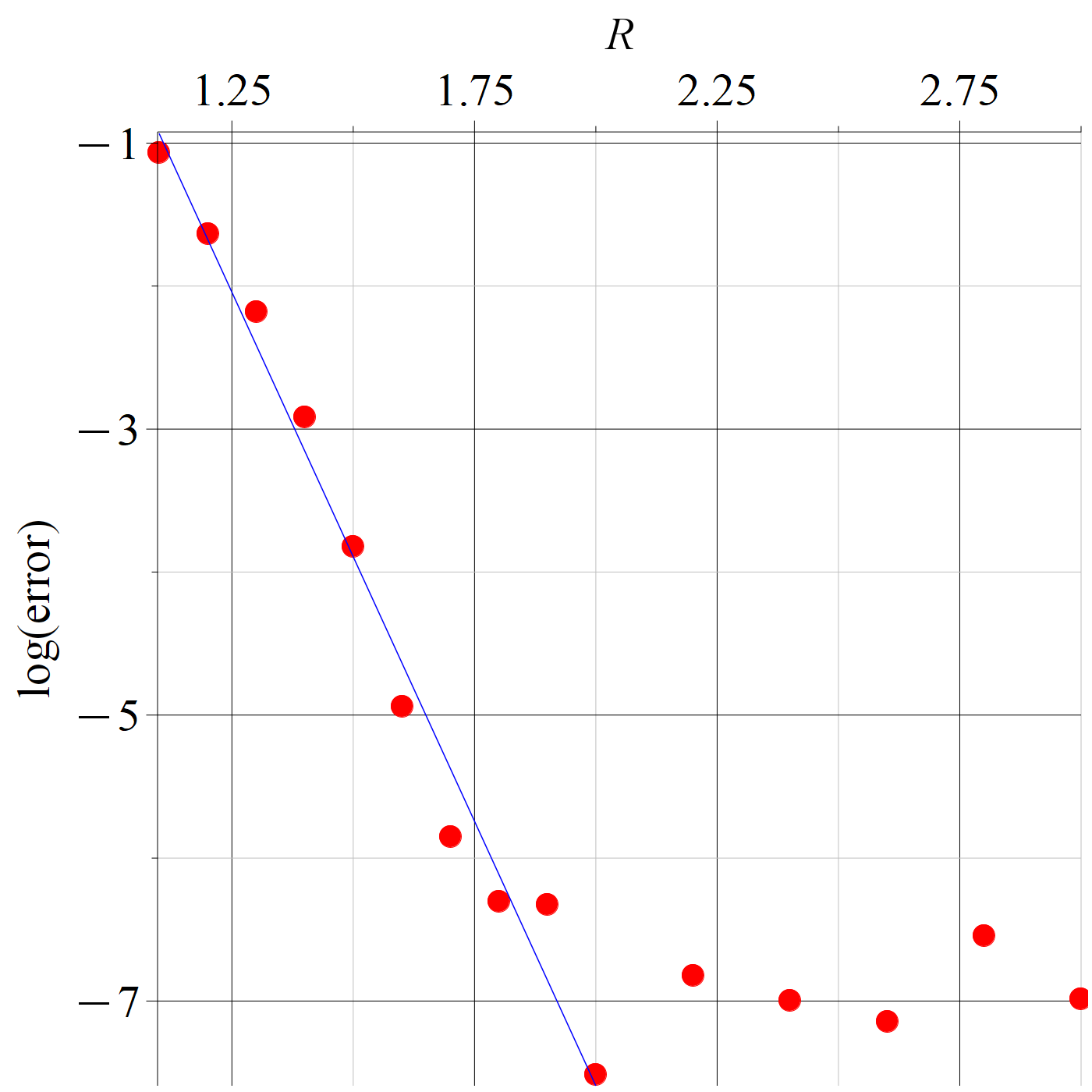}}
\caption{Left panel: supremum of error $|u(5,r,\theta)-u_{pml}(5,r,\theta)|$ for $0\leq r\leq 1$ as function of truncation
radius $R$. Right panel: Logarithmic plot of this error.}
\label{fig:pmlsol3}
\end{figure}

Figure \ref{fig:pmlsol3} illustrates the effect of varying the truncation radius $R$ of the disk on which the PML system is solved.
The left panel graphs the maximum relative error in $|u(5,r,\theta)-u_{pml}(5,r,\theta)|$ versus $R$, while the right panel shows this same information on a logarithmic scale.
Note that for $R$ up to about $R=1.9$ this error decays approximately exponentially (as evidenced by the straight line
behavior in the right panel), in accord with Remark \ref{rem:mainresult2}. For $R>1.9$ the error in the FEM solution itself dominates and increasing $R$ provides no improvement.

\subsection{Experiment 3}

\begin{figure}[h]
\centerline{\includegraphics[width=0.5\textwidth]{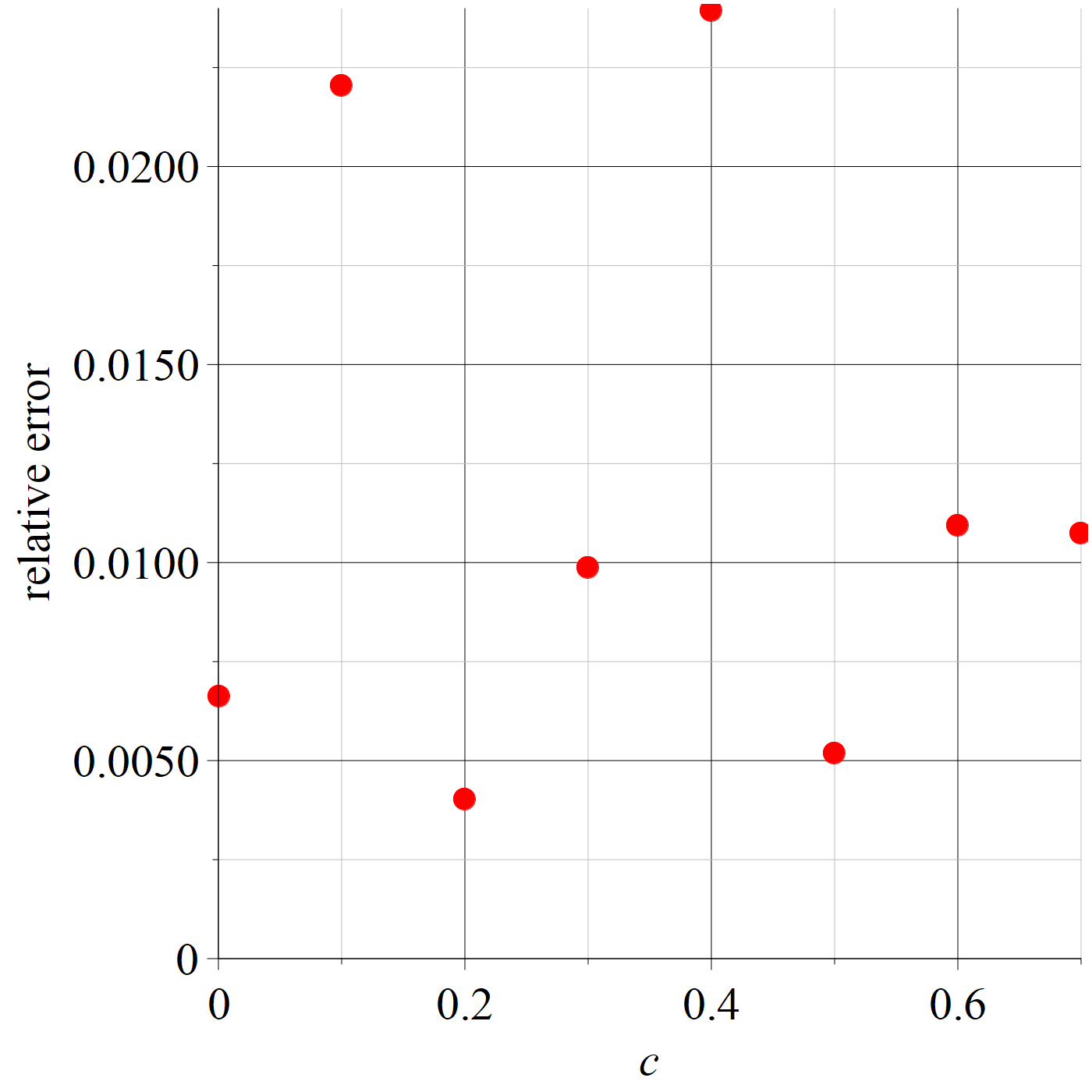}}
\caption{PML error as a function of offset $c$.}
\label{fig:offsetL2err}
\end{figure}

In this computation we use a source function
\[
f(t,r,\theta) = 1_{\{r'<0.25\}}
\]
but offset from the origin by taking $r'=\sqrt{(x-c)^2+y^2}$ for $c=0,0.1,0.2,\ldots,0.7$.
Figure \ref{fig:offsetL2err} shows the relative error in the PML solution, again as measured in the supremum norm,
with the truncation value $R=2$ in the PML equations.
We note that the analysis leading to Proposition \ref{prop:mainresult} and Remark \ref{rem:mainresult2} remains
valid so long as the support of $f$ in (\ref{wave-eq}) is contained in $r<1$, as it is here.
The error seems to have no dependence on the offset $c$, even as the support of $f$ (when $c=0.7$) comes quite
close to the boundary of the unit disk. This error is primarily due to the FEM discretization and first-order
scheme in time. 

%% file: Appendix1.tex
\section{Appendix}
\label{sec:appendix}

\subsection{Some symmetries}

Due to the integral representations for Bessel functions (see for instance \cite{Watson} section )it is fairly straightforward to verify that the functions $H_n^{(1)}(z)$ and $J_n(z)$, when considered for complex arguments, satisfy the following symmetry relations across the imaginary axis
\begin{lemma}
\label{relcompl}
Let $n$ be a non-negative integer. Then
$$
J_n(-\overline{z})= (-1)^n \overline{J_n(z)}~,~ \hbox{ for all } z\in \C
$$
and
$$
H_n^{(1)}(-\overline{z})= (-1)^{n+1} \overline{H_n^{(1)}(z)}~,~ \hbox{ for all }~z \in \C_+ =\{ z\in \C~ \hbox{with  Im}(z)>0 \}.
$$
\end{lemma}
\vskip 5pt
\noindent
Based on the symmetries asserted in Lemma \ref{relcompl}, simple calculations give the following lemma
\begin{lemma}
\label{relcompl2}
Let $n$ be a non-negative integer, and let $F_n$ and $G_n$ denote the functions
$$
F_n(z)= z\frac{d}{dz}J_n(z)~,~~z\in \C~,  ~~~\hbox{ and }~~~ G_n(z) = z\frac{d}{dz} H_n^{(1)}(z)~ ,~~z\in \C_+~.
$$
Then
$$
F_n(-\overline{z})= (-1)^n \overline{F_n(z)}~,~ \hbox{ for all } z\in \C~,
$$
and
$$
G_n(-\overline{z})= (-1)^{n+1} \overline{G_n(z)}~,~ \hbox{ for all }~z \in \C_+.
$$
\end{lemma}
\vskip 5pt
\noindent
Finally from Lemma \ref{relcompl} and the fact that $|J_n(z)|\le C_n (|z|+1)^{-1/2}$, $z\in \C_+$, we also obtain
\begin{lemma}
\label{relcompl3}
Let $n$ be a non-negative integer, then the function $g_n$ defined by
$$
g_n(z)=\int_0^1J_n(zs)s\,ds~~~z\in \C_+~,
$$
satisfies
$$
g_n(-\overline{z})= (-1)^n \overline{g(z)}~~~\hbox{ and }~~~ |g_n(z)|\le C_n (|z|+1)^{-1/2}~,~ z\in \C_+  ~.
$$
\end{lemma}
\vskip 5pt
\noindent
We leave the proofs of these three lemmata to the reader.

\subsection{Details of the numerics}

To solve (\ref{vpmleq})-(\ref{upmleq}) numerically we first convert these equations to a first-order system in time by introducing a function $q = \partial u_{pml}/\partial t$ and splitting (\ref{upmleq}) into a pair of equivalent equations
\begin{equation}
\label{eqn:splitupml}
\begin{aligned}
\frac{\partial q}{\partial t} & = -2\alpha q - \alpha^2 u_{pml} + \frac{1}{r^2}\frac{\partial^2 u_{pml}}{\partial \theta^2} + \frac{w_{pml}}{r^2} + f\\
\frac{\partial u_{pml}}{\partial t} & = q
\end{aligned}
\end{equation}
where the first equation in (\ref{eqn:splitupml}) has been divided by $r^2$. The two equations in (\ref{eqn:splitupml}) in conjunction with (\ref{vpmleq})-(\ref{wpmleq}) form a set of four first-order equations that can be marched out in time.

For numerical solution in FEniCSx the equations are converted to rectangular ($xy$) coordinates. The equations
of (\ref{eqn:splitupml}) become
\begin{align}
\frac{\partial q}{\partial t} & = -2\alpha q - \alpha^2 u_{pml} +
\frac{(\langle -y,x\rangle\cdot \nabla )^2 u_{pml}}{x^2+y^2} + \frac{w_{pml}}{x^2+y^2} + f\label{eqn:pmlrec1}\\
\frac{\partial u_{pml}}{\partial t} & = q.\label{eqn:pmlrec2}
\end{align}
It's worth noting that in the region $r<1$ the second equation of (\ref{eqn:splitupml}) becomes merely
$\frac{\partial q}{\partial t} = \Delta u_{pml} + f$ or
\begin{equation}
\label{eqn:pmlrec2b}
\frac{\partial q}{\partial t} = \frac{\partial^2 u_{pml}}{\partial x^2} + \frac{\partial^2 u_{pml}}{\partial y^2} + f,
\end{equation}
thus avoiding any singularity at the origin.
Equation (\ref{vpmleq}) can be expressed as
\begin{equation}
\label{eqn:vpmlrect}
\frac{\partial v_{pml}}{\partial t} = -\beta v_{pml} + \langle x,y\rangle \cdot \nabla q + \alpha \langle x,y\rangle \cdot \nabla u_{pml}
\end{equation}
by making use of $q=\frac{\partial u_{pml}}{\partial t}$ and $r\partial_r = \langle x,y\rangle\cdot \nabla$.
Finally, define an auxiliary function $\eta=\partial v_{pml}/\partial t$ but in the form
\begin{equation}
\label{eqn:etadef}
\eta = -\beta v_{pml} + \langle x,y\rangle \cdot \nabla q + \alpha \langle x,y\rangle \cdot \nabla u_{pml}
\end{equation}
so that $\frac{\partial^2 v_{pml}}{\partial r\partial t} = \frac{\partial\eta}{\partial r}$. Then (\ref{wpmleq})
can be expressed as
\begin{equation}
\label{eqn:wpmlrect}
\frac{\partial w_{pml}}{\partial t} = -\beta w_{pml} + \langle x,y\rangle \cdot \nabla\eta + \alpha \langle x,y\rangle \cdot \nabla v_{pml}.
\end{equation}
Equations (\ref{eqn:pmlrec1}), (\ref{eqn:pmlrec2}), (\ref{eqn:vpmlrect}), (\ref{eqn:etadef}), and (\ref{eqn:wpmlrect})
(replacing (\ref{eqn:pmlrec2}) with (\ref{eqn:pmlrec2b}) when $r<1$) govern the evolution of the system in time.

To solve the system numerically we form a weak version of these equations by multiplying each equation by a suitable test function $\psi(x,y)$ that vanishes at $r=R$ and integrating by parts. We also approximate each time derivative with a backward Euler difference, for example, $\frac{\partial q}{\partial t}(t_n) \approx (q^{n+1}-q^n)/\Delta t$, where $\Delta t$ is a fixed time step,
$t_n=n\Delta t$, and $q^n$ indicates an approximation to $q(t_n)$. With the notation
$(p,q)=\int_{D_R}p(x)q(x)\,dx$ ($D_R$ a disk of radius $R$ centered at the origin) the resulting system can be expressed as
\begin{align}
& (\left (1 + 2\alpha\,\Delta t) q^{n+1},\psi \right ) = (q^n,\psi) - (\alpha^2 u_{pml}^n,\psi )\,\Delta t
+ \left (\frac{w_{pml}^n}{x^2+y^2},\psi \right )\,\Delta t\nonumber\\
& - (\langle -y/r,x/r\rangle \cdot \nabla u_{pml}^n,\langle -y/r,x/r\rangle \cdot \nabla \psi)\,\Delta t + (f^{n+1/2},\psi)\,\Delta t\label{eqn:pml1rectweak}\\
& \nonumber \\
& (u_{pml}^{n+1},\psi) = (u_{pml}^n,\psi) + (q^{n+1},\psi)\,\Delta t\label{eqn:pml2rectweak}\\
& \nonumber \\
& ((1+\beta\,\Delta t)v_{pml}^{n+1},\psi) = (v_{pml}^n,\psi) + (\langle x,y\rangle \cdot \nabla q^{n+1},\psi)\,\Delta t + (\alpha \langle x,y\rangle \cdot \nabla u_{pml}^{n+1},\psi)\,\Delta t\label{eqn:pml3rectweak}\\
& \nonumber \\
& \eta^{n+1} = -\beta v_{pml}^{n+1} + \langle x,y\rangle \cdot \nabla q^{n+1} + \alpha \langle x,y\rangle \cdot \nabla u_{pml}^{n+1}\label{eqn:pmletarectweak}\\
& \nonumber \\
& ((1+\beta\,\Delta t)w_{pml}^{n+1},\psi) = (w_{pml}^n,\psi) - 2(\eta^{n+1},\psi)\,\Delta t - (\eta^{n+1},x\psi_x+y\psi_y)\,\Delta t\nonumber\\
& + (\alpha \langle x,y\rangle \cdot \nabla v_{pml}^{n+1},\psi)\,\Delta t.\label{eqn:pml4rectweak}
\end{align}
with $r=\sqrt{x^2+y^2}$; we also used a centered difference for $f$. For (\ref{eqn:pmlrec2b}) (used when $r<1$) we obtain
\begin{equation}
\label{eqn:eqnpml1rectweakb}
(q^{n+1},\psi) = (q^n,\psi) - (\nabla u_{pml}^n,\nabla \psi) + (f^{n+1/2}),\psi)\,\Delta t
\end{equation}

For each choice of $n=0,1,\ldots$ we form the relevant linear system for (\ref{eqn:pml1rectweak})-(\ref{eqn:eqnpml1rectweakb})
and solve for the updates $u_{pml}^{n+1},\ldots,w_{pml}^{n+1}$. These computations are performed
in the FEniCSx software with a basis of quadratic Lagrange elements, on a triangularization of size approximately $15000$.
The resulting system is then solved to form $q^{n+1}$.
Then (\ref{eqn:pml2rectweak}) is solved for $u_{pml}^{n+1}$,
(\ref{eqn:pml3rectweak}) is solved for $v_{pml}^{n+1}$, $\eta^{n+1}$ is computed using (\ref{eqn:pmletarectweak}), and
(\ref{eqn:pml4rectweak}) is solved for $w_{pml}^{n+1}$, all with the same parameters.
We solve out to time $T=5$ using $1000$ steps of size $\Delta t = 0.005$.
%\marginpar{There is a von Neumann stability analysis that could be performed here.} 